\documentclass{amsart}
\newtheorem{theorem}{Theorem}[section]

\newtheorem{proposition}[theorem]{Proposition}
\newtheorem{corollary}[theorem]{Corollary}
\theoremstyle{definition}

\theoremstyle{remark}
\newtheorem{remark}[theorem]{Remark}

\numberwithin{equation}{section}

\newcommand{\loc}   {loc}

\newcommand{\ERRE}{\mathbb{R}}
\newcommand{\ENNE}{\mathbb{N}}
\newcommand{\Z}{\mathbb{Z}}

\newcommand{\dive}{\mathop{div}}
\newcommand{\F}{\mathbb{F}}
\newcommand{\B}{\mathbb{B}}
\newcommand{\A}{\mathbb{A}}
\newcommand{\U}{\mathbb{U}}
\begin{document}   
\title
[Hyperbolic to Parabolic Relaxation]
{Convergence
of Singular Limits \\ for Multi-D Semilinear
Hyperbolic Systems\\ to Parabolic Systems  }
\author{Donatella Donatelli}
        \address{Dip. di Matematica Pura ed Applicata\\
       Universit\`a degli Studi dell'Aquila\\
         67100 L'Aquila, Italy}
        \email{donatell@univaq.it}
\author{Pierangelo Marcati}
         \address{Dip. di Matematica Pura ed Applicata\\
       Universit\`a degli Studi dell'Aquila\\
       67100 L'Aquila, Italy}
        \email{marcati@univaq.it}

\begin{abstract}
In this paper we investigate the zero-relaxation limit of the
following multi-D semilinear hyperbolic
system in pseudodifferential form:
$$W_{t}(x,t) + \frac{1}{\varepsilon}A(x,D)W(x,t)=
\frac{1}{\varepsilon ^2} B(x,W(x,t))+\frac{1}{\varepsilon}
D(W(x,t))+E(W(x,t)).$$
We analyse the singular convergence, as $\varepsilon \downarrow 0$, in the 
case which leads to a limit system of parabolic type.
The analysis is carried out by using the following steps:
\begin{itemize}
\item[(i)] We single out algebraic ``structure conditions'' on the full system,
motivated by formal asymptotics, by some examples of discrete velocity models
in kinetic theories.
\item[(ii)] We deduce ``energy estimates '', uniformly in $\varepsilon$,
by assuming the existence of a symmetrizer having the so called block
structure and by assuming ``dissipativity conditions '' on $B$.
\item[(iii)] We perform the convergence analysis by using generalizations
of Compensated Compactness due to Tartar and G\'erard.
\end{itemize}
Finally we include examples which show how to use our theory to 
approximate prescribed general quasilinear parabolic systems, satisfying 
Petrowski parabolicity condition, or general reaction diffusion systems.
\medbreak
\medbreak
\end{abstract}
\maketitle
\begin{section}{Introduction}
In this paper we study the following semilinear multidimensional
hyperbolic system with  a small parameter $\varepsilon >0$
{\small
\begin{equation}
W_{t}(x,t) \!+ \!\frac{1}{\varepsilon}A(x,D)W(x,t)=
\frac{1}{\varepsilon ^2} B(x,W(x,t))+\frac{1}{\varepsilon}
D(W(x,t))+E(W(x,t)), \label{0}
\end{equation}}
where $W=W(x,t)$ belongs to $\ERRE ^N$, $x\in \ERRE ^d$,
$t\geq 0$, $A(x,D)$ is a first order pseudodifferential operator. The system 
(\ref{0}) includes also the case of  hyperbolic differential
operators of the form
{\small
\begin{equation}
W_{t}(x,t) + \frac{1}{\varepsilon}\sum_{j =
1}^{d}A_{j}(x)\partial_{j}W(x,t)= \frac{1}{\varepsilon ^2}
B(x,W(x,t))+\frac{1}{\varepsilon} D(W(x,t))+E(W(x,t)), \label{00}
\end{equation}}
where $A_{j}(x)$, $j=1, \ldots,d $ are $N\times N$ matrices for
any $x \in \ERRE ^d$. Our aim is to describe  the limiting
behaviour of the system (\ref{0}) as $\varepsilon$ goes to zero.
We look for structure condition to ensure that (\ref{0})
approximate a second order parabolic system. Our  interest in this
problem is  motivated also by  a very strong similarity with the
limiting structure appearing in the investigation of the
hydrodynamic limit for the  Boltzmann equation, in particular in
the discrete velocity case. The Boltzmann equation describes the
evolution of the density $f(x, \xi, t)$ of particles which are at
time $t$ in position $x$ with velocity $\xi$ and has the following
form $$\nu f_t +\xi \nabla \cdot f = \frac{1}{\varepsilon}Q(f,f)$$
where $\nu$ is the Mach number and  $\varepsilon$ the Knudsen
number. By averaging $f(x, \xi, t)$ in $\xi$ and by using higher
order momenta we can define a hierarchy of macroscopic quantities.
The investigation of the hydrodynamic limit regards the behaviour
of those quantities (actually combined with the closure problem)
as the Knudsen number goes to zero. In the case where the Mach
number is of the same order of the Knudsen number our limit can be
described by the Navier Stokes equation, otherwise when  the Mach
number is fixed and the Knudsen number tends to zero we end up
with the Euler equation. Those limits allow us to understand the
differences between relaxation limits of hyperbolic type to
parabolic with respect to those one of hyperbolic to hyperbolic
type \cite{CLL94}, \cite{KT96}, \cite{KT97}, \cite{KT99},
\cite{Liu87}, \cite{Yo99}, \cite{GoRe}. In particular if we deal with a
discrete velocity models, the equivalent `` Boltzmann equation''
is a semilinear hyperbolic system. The simplest example is given
by the Carleman's equation $$ \left\{ \begin{array}{ll}
f_{1t}+\displaystyle{\frac{1}{\varepsilon}}f_{1x}=
\displaystyle{\frac{1}{\varepsilon ^2}}(f^2_2-f^2_1)\\ \\
f_{2t}-\displaystyle{\frac{1}{\varepsilon}}f_{2x}=
\displaystyle{\frac{1}{\varepsilon ^2}}(f^2_1-f^2_2),
\end{array}
\right. $$ where we take $\xi \in\{-1,1 \}$ and $f_1=f(x,1,t)$,
$f_2=f(x,-1,t)$. By rescaling the variable we get that $\rho =f_1 +f_2$
as $\varepsilon \downarrow 0$ satisfies the nonlinear diffusion
equation $$\rho _t=\frac{1}{2}(\log \rho)_{xx}.$$ This asymptotic problem was
first investigated by  Kurtz \cite{Kur73} and McKean \cite{Mck75}.
Therefore the nonlinear diffusion problem, obtained as the limit
of the Cattaneo hyperbolic nonlinear heat conduction equation, was
proved by Marcati, Milani and Secchi \cite{MMS88}. The paper of
Marcati and Milani \cite{MM90} concernes the pourous media flow as
the limit of the Euler equation in 1-D, later generalized by
Marcati and Rubino \cite{MR} to the multi-D case. Relaxation
phenomena of the same nature appear in the zero relaxation limits
for the Euler-Poisson model for the semiconductor devices and it
was investigated by Marcati and Natalini \cite{MN95a},
\cite{MN95b} in the 1-D case and by Lattanzio and Marcati
\cite{LaM99} in the multi-D case. More recently Lions and Toscani
\cite{LT96} investigated a discrete velocity model leading to the
pourous media flow. All of these papers, with exception of
\cite{Kur73}, \cite{Mck75}, make use of the techniques of
compensated compactness. Similar ideas have been applied by
Marcati and Rubino \cite{MR} to show the general theory for $2
\times 2$ systems in the 1-D case and to propose a general
framework that we are going to investigate here in the semilinear
system case. Models of BGK type approximation have been
successfully studied in this framework by Bouchut, Guarguaglini and Natalini
\cite{BGN} and Lattanzio, Natalini \cite{LaN02} in the case of 1-D
systems. The framework of \cite{MR} was also investigated in the
quasilinear case by Lattanzio and Yong \cite{LaY} for $H^s$-
smooth solutions. 
Recently, with a similar approach, Junk and Yong \cite{JY} derived the incompressible 
Navier Stokes equations form the BGK model.
Preliminary results concerning semilinear
systems have been obtained  in \cite{MD99}, \cite{MD00a}, in
particular in \cite{MD00a} we considered  a 1- D semilinear system
with variable coefficients. Already in that case the classical
compensated compactness is not sufficient and it is necessary to
use a generalization of this theory due to Tartar \cite{TarM} and
P.G\'erard \cite{Ger}. From the technical point of view the
problem, here, has additional complications , since we deal with 
a  multi- D pseudodifferential system and since the  symmetrizers
are pseudodifferential operators.\\ The plan of the paper is as
follows. In Section 2 we give the basic definitions, the basic
notations and we recall some mathematical tools needed in this
paper. In Section 3 we describe the structure condition on our system
and on its symmetrizers and we will describe the formal framework of the
limiting process. The Section 4 is
devoted to give rigorous proof of the previous formal analysis,
the assumptions of the previous  section will allow us to obtain 
the a priori estimates of energy type uniformly in $\varepsilon $. 
By using only
the informations coming from the energy estimates and by means of the 
previous mentioned compactness
framework we will be  able to obtain  our relaxation results. The
structure conditions used in the limit process allow us also to
satisfy the parabolicity condition for the limit system in the
case of system (\ref{00}) . In Section 5 we show how this theory
can be greatly simplified in the case of constant coefficients.
Finally, in Section 6 we show how it is possible to use the theory
developed in the previous sections in order to approximate a given
parabolic system by means of a larger hyperbolic system with simpler 
nonlinear structure. We are able to approximate the following 
nonlinear parabolic systems in divergence form 
$$
U_t +\sum _{i=1}^{d}\partial _{i} \left( F_{i}(U)-\sum _{j=1}^{d}B_{ij}(U)\partial 
_{j}U\right )=G(U) 
$$
where $x \in \ERRE ^d$, $t \in \ERRE_+$, $U=U(x,t) \in \ERRE ^k$.
With the previous techniques we approximate also reaction- diffusion systems 
of the form
$$
U_t=\sum_{j,k=1}^{d}A_{j,k}(x)\partial _j\partial _k U +f(U)
$$
where $x \in \ERRE ^d$, $t \in \ERRE_+$, $U=U(x,t) \in \ERRE ^k$. 
In this latter case we use two different approach. The former is based on  the theory of 
pseudodifferential operators while the latter is based on symmetric
differential operators and 
can be more usefull for numerical computations. 
Let us remark that this latter approximation extends in multi-D those 
proposed by Jin and Liu \cite{JL98} and Lattanzio and Natalini \cite{LaN02} 
(example 5.3). Our theory does not require the use of $L^{\infty}$
estimates (for instance via invariant domains like Serre \cite{Ser00}). 
\end{section}

\begin{section}{Preliminary notions }

We start this section  introducing  the main notations and definitions
used in this article. In particular   we recall some basic facts and 
notations concerning  the
theory of pseudodifferential operators and we also recall our main
compactness tools used in the  strong convergence analysis. Namely
\begin{itemize}
\item[\bf(a)] $(\cdot, \cdot)$ denotes the scalar product in $\ERRE^q$,
$(q=1,2,...)$ and $| \cdot|$ the usual norm of $\ERRE^q$
$(q=1,2,...)$,
\item[\bf(b)] $\mathcal{M}_{m\times n}$ denotes the linear space of $m\times n$
matrices,
\item[\bf(c)] $( \cdot, \cdot )_2 $ denotes the scalar product in
$L^{2} (\ERRE^d )$ and $\|\cdot\|$ the norm in $L^{2}
(\ERRE^d)$,
\item[\bf(d)] $\mathcal{D}(\ERRE ^d \times \ERRE_+)$
denotes the space of test function
$C^{\infty}_{0}(\ERRE^d \times \ERRE_+)$, $\mathcal{D}'(\ERRE^d \times
\ERRE_+)$ the Schwarz space of distributions and $\langle \cdot, \cdot \rangle$
the duality bracket in $\mathcal{D}'(\ERRE^d \times \ERRE_+)$,
\item[\bf(e)]  $H$ is a separable Hilbert space, $ \mathcal{L}(H)$ the space of
bounded operators,
$ \mathcal{K}(H)$ the space of compact operators,
\item[\bf(f)] we denote by $H^{s}_{loc}(\Omega, H)$ the classical local Sobolev
space of order $s$, i.e.
$u\in H^{s}_{loc}(\Omega, H) \Longleftrightarrow \forall\: \varphi
\in C^{\infty}_{0}, \: (\widehat{\varphi u})\in L^2(\ERRE^n,
(1+|\xi|^2)^s d\xi)$.
\item[\bf(g)] we denote by $F_\nu \!(\!Z^{I},Z^{II}\!)$ the derivative respect
to the variable $Z^{II}$ \!of\! $F(\!Z^{I},\!Z^{II}\!)$.
\end{itemize}

We shall also make use of the notion of parabolicity for systems of
equations in various way (see Taylor \cite{Tay96} volume III,
\cite{Tay91}, Eidel'man \cite{Eid69}, Kreiss and Lorenz
\cite{KL89}). Let us consider the following system
\begin{equation}
u_t =\sum_{j,k} A^{j,k}(t,x,D'_x u)\partial _j \partial _k u
+B(t,x,D'_x u), \label{s}
\end{equation}
where $u\in \ERRE^p$, $A^{j,k}(t,x,D'_x u) \in \mathcal{M}_{p\times
p}$, $B(t,x,D'_x u)\in \ERRE^p$ and $D' _x$ is a differential operator
of order not greater than two. The system is said {\em strongly
parabolic} if there exists $c_0>0$, such that for all $\xi \in \ERRE^d$
one has $$ \sum_{j,k} A^{j,k}(t,x,D'_x u)\xi _j \xi _k \geq c_0
|\xi|^2 I.$$ Namely, if we denote $\displaystyle L(t,x,D'_x u,
\xi)=-\sum_{j,k} A^{j,k}(t,x,D'_x u)\xi _j \xi _k $ this condition is
equivalent to say $L+L^T$ is a negative definite matrix. Unfortunately,
this condition is often difficult to be verified, then we recall 
a more general notion often referred as {\em Petrowski
parabolicity} (see Taylor \cite{Tay96}
volume III, \cite{Tay91}).\\ We say that the system (\ref{s}) is {\em parabolic}
if, denoted by $\lambda _k(t,x,D'_x u, \xi)$ the eigenvalues of the
matrix $L(t,x,D'_x u, \xi)$, one has there exists $\alpha _0 >0$ such
that, for all $\xi \in \ERRE^d$, $$Re \lambda _k(t,x,D'_x u, \xi) \leq
-\alpha _0|\xi|^2.$$ The latter notion of parabolicity is equivalent
to ask the existence of a symmetric matrix $P_0(t,x,D'_x u, \xi)$, positive
definite (i.e. $P_0\geq cI>0$), homogeneous of degree $0$ in $\xi$,
such that $$-(P_0L+L^{\ast}P_0) \geq c |\xi|^2 I>0.$$ Let us recall
the basic notations concerning pseudodifferential operators to be used
later on, we refer to \cite{Tay91} for details. Assuming $\rho, \delta
\in [0,1],\ m \in \ERRE$, we denote $S^{m}_{\rho, \delta}$ the set of
$C^\infty$ symbols satisfying $$\left| D^\beta _x D^\alpha _\xi
p(x,\xi)\right|\leq C_{\alpha, \beta}\langle \xi \rangle
^{m-\rho|\alpha|+\delta|\beta|}$$ for all $\alpha, \beta$, where
$\langle \xi \rangle =(1+|\xi|^2)^{1/2}$. In such case we say that the
associated operator defined by $$P(x,D)f(x)=\int p(x, \xi) \hat f
(\xi) e^{i x \xi}d\xi:=OP(p(x,\xi))$$ (where$\hat f
(\xi)=(2\pi)^{-n}\int f(x) e^{-ix \xi}dx $ denotes the Fourier transform of
the function $f$) belongs to $OPS^{m}_{\rho, \delta}$. If there
are smooth symbols $p_{m-j}(x,\xi)$, homogeneous in $\xi$ of degree $m-j$ for
$|\xi|\geq 1$, i.e.  $p_{m-j}(x,r\xi)=r^{m-j}p_{m-j}(x,\xi)$ for $r,\
|\xi|\geq 1$, and if $$p(x,\xi)\sim \sum_{j\geq 0} p_{m-j}(x,\xi)$$ in
the sense that $$p(x,\xi)-\sum_{j\geq 0}^{N} p_{m-j}(x,\xi)\in
S^{m-N}_{1,0}
$$ for all $N$, then we say $p(x,\xi)\in S^m$. We define also
$$S^{-\infty}=\bigcap _{m>0}S^{-m}_{\rho, \delta}.$$ The following
properties will also be used here.
\begin{theorem}{\rm (\cite{Tay91}, Theorem 0.5 A)}
If $P(x,D) \in OPS^{0}_{\rho, \delta} $ and $0\leq \delta <\rho \leq
1$ then
$$P(x,D) : L^2(\ERRE^n) \longrightarrow L^2(\ERRE^n). $$
\label{p1}
\end{theorem}
\begin{theorem}{\rm (\cite{Tay81}, \!Chapter II, \!Lemma 6.2 Ex. \!\!\!8.1,
Chapter III, Proposition 3.1)}\\If $P(x,D) \in OPS^{0}_{\rho, \delta}
$, $\delta <\rho$ and if $\Re e p(x,\xi) \geq c>0$ then there  exists
$\widetilde{R}(x,D) \in OPS^{0}_{\rho, \delta}$ such that
$\widetilde{R}(x,D) \geq \eta I >0$, with $\Re e
P(x,D)=\frac{1}{2}(P+P^\ast)$,
 $$\Re e P(x,D) -\widetilde{R}(x,D) \in OPS^{-\infty}.$$ \label{p2}
\end{theorem}
\begin{proposition}{\bf (Commutator estimate)}
{\rm(\cite{Tay91}, Proposition 3.6 B)}\\
Given $P(x,D) \in OPS^{0}_{1,0} $ we have $$\|[P,\
f]u\|_{W^{\sigma,p}}\leq c \|f\|_{Lip}\|u\|_{W^{\sigma -1,p}} \qquad
\text {for $0\leq \sigma \leq 1$},$$ in particular if $\sigma =0$ and
$p=2$ $$\|[P,\ f]u\|_ {L^2(\ERRE^n)}\leq c \|f\|_{Lip}\|u\|_{H^{
-1}}. $$ \label{p3}
\end{proposition}
Finally we state here our main tool, due to Tartar and G\'erard
(\cite{TarM}, \cite{Ger})), to study the convergence of
quadratic forms with variable coefficients. The classical results concerning the use
of Compensated Compactness in the theory of hyperbolic systems are
reported in the books of Dafermos \cite{Daf00} and Serre \cite{Ser96}.
Let $H$, $H^\sharp$ denote separable Hilbert spaces, $\Omega \in \ERRE
^n$, an open set, $m\in \ENNE$. We have the following theorem taken
from G\'erard \cite{Ger}.
\begin{theorem}
{\bf (Compensated Compactness)}\\ Let $P\in OPS^m$ with principal
symbol $p(x,\xi)$ and $\{u_k\}$ be a bounded sequence of
$L^{2}_{loc}(\Omega, H)$, such that $u_k \rightharpoonup u $.  Assume
that there exists a dense subset $D\in H^\sharp$ such that, for any
$h\in D$, the sequence $(\langle Pu_k,h \rangle )$ is relatively
compact in $H^{-m}_{loc}(\Omega)$. Moreover, let $q\in C(\Omega,
\mathcal{K}(H))$.
\begin{itemize}
\item[{\bf (i)}] If $q=q^\ast$ and for all $(x,\xi,h)\in
\Omega \times S^{n-1} \times H$, one has
$$
(p(x,\xi)h=0)\;\Rightarrow \; (\langle q(x)h,h\rangle \geq 0)$$ then,
for any nonnegative $\varphi \in C^{\infty}_{0}(\Omega)$
$$\liminf_{k\rightarrow \infty}\int_{\Omega}\varphi\langle
q(x)u_k,u_k\rangle dx \geq \int_{\Omega}\varphi \langle q(x)u,u\rangle
dx.$$ 
\item[{\bf (ii)}] If for all $(x,\xi,h)\in
\Omega \times S^{n-1} \times H$, one has 
$$(p(x,\xi)h=0)\;\Rightarrow \; (\langle q(x)h,h\rangle =0)$$
then 
$$\langle q(x)u_k,u_k\rangle \quad \text{ converges to } \quad
\langle q(x)u,u\rangle \qquad \text{in $\mathcal{D}'(\Omega)$}.$$
\end{itemize}
\label{g}
\end{theorem}
The previous theorem holds also if instead of a differential operator
$P\in OPS^m$ we consider the differential operator $$
Pu(x)=\sum_{|\alpha|\leq m}\partial ^\alpha ( a_\alpha (x)u(x)),
$$ where $\alpha \in \ENNE^n$, $|\alpha|\leq m$, $a_\alpha \in
C(\Omega , \mathcal{L}(H,H^\sharp))$.
\end{section}
\begin{section}{Multidimensional Framework}
We shall restrict our analysis to the case $E(W)\equiv  0$ to simplify the 
computations, but our results easily extend to the case of 
$E(\cdot)\in Lip(\ERRE^{N},\ERRE^{N})$.   
\begin{subsection}{Decoupled system}

We will consider the following semilinear system of equations
\begin{equation}
W_{t}(x,t) + A(x,D)W(x,t) = B(x,W(x,t))+D(W(x,t)) \label{1}
\end{equation}
where $t\geq 0$, $x\in \ERRE ^d$, $W\in \ERRE ^N$. We assume
the following hypotheses hold.
\begin{itemize}
\item[{\bf(A.1)}] $B(\cdot,\cdot)\in C^{1}(\ERRE ^{d+N}, \ERRE ^N)$,
$D(\cdot)\in C^{1}(\ERRE ^{N}, \ERRE ^N)$
\item[\bf(A.2)] $A(x,D)\in OPS^1$, the system (\ref{1}) is hyperbolic,
namely the principal symbol of $A(x,D)$ is the matrix $ A(x,\xi) \in
\mathcal{M}_{N\times N}$ whose eigenvalues for $(x,\xi) \in
\mathbb{R}^{d}\times \ERRE ^d$, $\xi \neq 0$ are pure imaginary,
\item[\bf(A.3)] denote by $S=span \left\{B(x,W) ~|~ (x,W)\in \ERRE
^{d+N}\right\}$  then $dim\ S =N-k$, $0<k<N$.
\end{itemize}
\begin{remark}
We point out that (\ref{1}) includes the case of the following
hyperbolic semilinear system $$ W_{t}(x,t) + \sum_{j =
1}^{d}A_{j}(x)\partial_{j}W(x,t) = B(x,W(x,t))+D(W(x,t))$$ where
$A_{j}(x)\in \mathcal{M}_{N\times N}$, $j = 0, \ldots, d$ and for all
nonzero vector $\xi \in \mathbb{R}^{d}$, the matrix $ A(x,\xi) =
\displaystyle{\sum_{j = 1}^{d}\xi _{j}A_{j}(x)}$ has real
eigenvalues. 
\end{remark}
The aim of this section is to decouple the full system in order
to single out the conserved quantities from the others.
We want to motivate the idea of decoupling by considering Carleman's 
system:
\begin{equation}
 \left\{ \begin{array}{ll}
f_{1t}+\displaystyle{\frac{1}{\varepsilon}}f_{1x}=
\displaystyle{\frac{1}{\varepsilon ^2}}(f_2+f_1)(f_2-f_1)\\ \\
f_{2t}-\displaystyle{\frac{1}{\varepsilon}}f_{2x}=
\displaystyle{\frac{1}{\varepsilon ^2}}(f_2+f_1)(f_1-f_2).
\end{array}
\right. 
\label{a1}
\end{equation}
Setting $\rho =f_1+f_2$ and $m=\displaystyle{\frac{f_1-f_2}{\varepsilon}}$
we get $\rho, \ m$ satisfy 
\begin{equation}
 \left\{ \begin{array}{ll}
\rho _t +m_x=0\\ \\
\varepsilon ^2 m_t +\rho _x =-2\rho m.
\end{array}
\right. 
\label{a2}
\end{equation}
In this way we  decoupled (\ref{a1}) in   (\ref{a2}) and we 
have isolated $\rho$ that can be easily seen, to be  the conserved quantity.
Now we want to do the same procedure on our system (\ref{1}). 
Let us consider the hypothesis
(A.3), then there exists a matrix $P^I \in \mathcal{M}_{k\times N}$ such 
that $P^I B(W)=0 $ and $Z^I=P^I W$ is the conserved vector. Therefore,  
we can construct an
invertible matrix $P=\begin{bmatrix} P^I \\ P^{II}
\end{bmatrix}$,  $P^{II} \in \mathcal{M}_{(N-k)\times N}$. 
Now, for any $W\in \ERRE ^N$, we set
\begin{align*}
Z^I&= P^I W, &\hspace{-5cm}Z^{II}&= P^{II}W,\\ Z&=\begin{bmatrix} Z^I
    \\ Z^{II}
\end{bmatrix},
           &Q(Z^I,Z^{II})&=P^{II}B(P^{-1}Z)\\  
             D^I(Z^I,Z^{II})&=P^{I}D(P^{-1}Z) &
           \hspace{-5cm} D^{II}(Z^I,Z^{II})&=P^{II}D(P^{-1}Z)\\
            PA(x,\xi)P^{-1}&=
               \begin{bmatrix} M^{11}(x,\xi) & M^{12}(x,\xi)\\ 
               M^{21}(x,\xi) & M^{22}(x,\xi)\\ \end{bmatrix}=
              M(x,\xi) & 
              M^{ij}(x,D)&=OPM^{ij}(x,\xi),
\end{align*}
hence by using the previous notations 
we can rewrite the system (\ref{1}) in the following decoupled form
\begin{equation}
\left\{ \begin{array}{ll} Z^{I}_{t}+M^{11}(x,D) Z^{I}+
  M^{12}(x,D)Z^{II} = D^I(Z^I, Z^{II})\\ \\ Z^{II}_{t}+M^{21}(x,D)
  Z^{I}+ M^{22}(x,D)Z^{II}= Q(x,Z^I, Z^{II})+D^{II}(Z^I, Z^{II}),
  \end{array} \right.
\label{2}
\end{equation}
where by construction $Z^I=Z^I(x,t) \!\in\! \ERRE ^k$, $Z^{II}=Z^{II}(x,t) \in
\ERRE ^{N-k}$ and $M^{11}(x,\xi)\in\mathcal{M}_{k\times k}$,
$M^{12}(x,\xi)\!\in\! \mathcal{M}_{k\times(N-k)}$, $M^{21}(x,\xi)\!\in\!
\mathcal{M}_{(N-k)\times k}$,
$M^{22}(x,\xi)\!\in\!\mathcal{M}_{(N-k)\times(N-k)}$, $D^I(Z^I,Z^{II})\in
\ERRE ^k$, $D^{II}(Z^I,Z^{II})\in \ERRE ^{N-k}, Q(x,Z^I,Z^{II})\in
\ERRE ^{N-k}$. The previous transformation does not affect the hyperbolic 
character of our system.
\end{subsection}
\begin{subsection}{Structural conditions}
In order to perform our analysis on system (\ref{1}) we need the
following structural assumption:
\begin{itemize}
\item[\bf(S.1)] $M^{11}(x,\xi)=0$, for any $(x,\xi)\in
\ERRE^d\times\ERRE^d$.
\end{itemize}
This condition is natural if we consider system (\ref{a2}). In that case
$M(x,\xi)$ is antisymmetric and is given by
$$M(x,\xi)=\begin{pmatrix}
                      0 & 1\\
                      1 & 0 
                     \end{pmatrix}. $$
The same hypothesis is assumed also by  Lions and Toscani
\cite{LT96}. It is also implicitly  contained in the work of 
Marcati and Rubino \cite{MR}. In fact in  \cite{MR} the relaxation of the 
following quasilinear nonhomogeneous $2\times 2$ hyperbolic system is 
considered,
\begin{equation}
	\begin{cases}
	w _{s } + f\left(w , z\right)_y = 0 & \\			
        z_{s } + g\left(w , z\right)_y = h\left(w , z\right), & 
        \end{cases}		
\end{equation}
where $y\in\mathbb{R}$, $s \geq 0$, with the assumption 
$$f(w, 0) = 0,$$
which in the linear case is equivalent to (S.1).
An analogous condition is assumed by Lattanzio and Yong in \cite{LaY}, where
they study the singular limits for the initial value 
problem
$$W_t+\sum_{j=1}^d \overline A_j(W) W_{x_j} +
\frac{1}{\varepsilon}\sum_{j=1}^d A_j(\varepsilon W) W_{x_j}=
\frac{Q(W)}{\varepsilon ^2},  $$
$$W(x,0)=W_0(x;\varepsilon),$$
by validating the formal asymptotic approximation in the framework of 
$H^s$- smooth solutions. In \cite{LaY}
they set
\begin{align*}
 \overline A_j=\begin{pmatrix}
                \overline A_j ^{11} & \overline A_j ^{12}\\
                 \overline A_j ^{21} & \overline A_j ^{22}   
               \end{pmatrix} & \qquad 
A_j=\begin{pmatrix}
                A_j ^{11} & A_j ^{12}\\
                A_j ^{21} & A_j ^{22}   
               \end{pmatrix} 
\end{align*} 
and the key structure condition is given by
$$ A_j ^{11}(0)=0, \qquad \text{for all $j$}.$$
In practice the condition (S.1) is essential otherwise  
the only relaxation process  would be the trivial one, which relaxes on
the null solution. In the case (S.1) is not valid, we actually, do not
have relaxation from hyperbolic to parabolic systems but we have a 
multiscale phenomena which involves the simultaneous action of distinct
relaxation mechanisms. To clarify this issue, we  follow 
the formal asymptotic approximations of  Lattanzio and Yong \cite{LaY}, 
namely we consider the scaled  system
\begin{equation}
W_t+
\frac{1}{\varepsilon}\sum_{j=1}^d A_j(x) W_{x_j}=
\frac{Q(W)}{\varepsilon ^2},
\label{a3}  
\end{equation}
where $W\in \ERRE$, $(x,t)\in \ERRE ^d \times \ERRE _+$, 
$A_j(x)$ are smooth $N\times N$ matrices with the following assumption
\begin{itemize}
\item[\bf (a)] the first $k$ components of $Q(W)$ are zero that is
$$Q(W)=\begin{pmatrix}
0\\
q(W)
\end{pmatrix}$$
with $q(W)\in \ERRE^{N-k}$.
\end{itemize}
Corresponding to this decompostion of $Q$, we set 
\begin{align*}
 W=\begin{pmatrix}
               u \\
               v   
               \end{pmatrix} & \qquad 
A_j=\begin{pmatrix}
                A_j ^{11} & A_j ^{12}\\
                A_j ^{21} & A_j ^{22}   
               \end{pmatrix}.  
\end{align*}
Furthermore we assume:
\begin{itemize}
\item[\bf (b)] $q(u,v)=0$ if and only if $v=0$, $q_v(u,0)$ is invertible
for any $u$.
\end{itemize}
We look for a solution of the form
$$O_\varepsilon (x,t) \sim \sum_{k=0}^{\infty} \varepsilon ^k O_k(x,t)=
O_0(x,t) +\sum_{k=1}^{\infty} \varepsilon ^k O_k(x,t).$$
Recalling (\ref{a3})  we introduce
$$R(W)=W_t+
\frac{1}{\varepsilon}\sum_{j=1}^d A_j(x) W_{x_j}-
\frac{Q(W)}{\varepsilon ^2} $$
Since $O_\varepsilon$ is expected to solve (\ref{a3}) we want 
$R(O_\varepsilon)=0$. 
Thus it follows 
\begin{align*}
Q(O_0)&=0\\
\sum_{j=1}^d A_j(x) O_{0x_{j}} &=Q_W(O_0)O_1\\
O_{kt} +
\sum_{j=1}^d A_j(x) O_{k+1,x_{j}}&= Q_W(O_0)O_{k+2}
\end{align*}
Set $\displaystyle{O_k =\begin{pmatrix} u_k \\v_k \end{pmatrix}}$. The above
equations can be rewritten as
\begin{align}
&q(u_0 ,v_0)=0\label{a5}\\
&\sum_{j=1}^d \left(A_j^{11}(x) u_{0x_{j}} +A_j^{12}(x) v_{0x_{j}}\right)
=0\label{a6}\\
&\sum_{j=1}^d \left(A_j^{21}(x) u_{0x_{j}} +A_j^{22}(x) v_{0x_{j}}\right)=
q_u(u_0,v_0)u_1 +q_v(u_0,v_0)v_1\label{a7}\\
&u_{kt}+
\sum_{j=1}^d \left(A_j^{11}(x) u_{k+1,x_{j}}+A_j^{12}(x) v_{k+1,x_{j}}\right)
=0\label{a8}\\
&v_{kt} +
\sum_{j=1}^d\left( A_j^{21}(x) u_{k+1,x_{j}}+A_j^{22}(x) v_{k+1,x_{j}}\right)
=
q_u(u_0,v_0)u_{k+2} \\ \notag &+q_v(u_0,v_0)v_{k+2}.\label{a9}
\end{align}    
From (\ref{a5}) and (b) it follows $v_0=0$, then (\ref{a6}) reduces to
$$\sum_{j=1}^d A_j^{11}(x) u_{0x_{j}}=0,$$
so $u_{0x_{j}}=0$. By using (b) and (\ref{a7}) we get $v_1=0$ and 
setting $k=0$ in (\ref{a8}) entails $u_1=0$. Up to now we found 
$u_0=v_0=u_1=v_1=0$. By assuming inductively $u_p=v_p=0$, by using 
the previous relations and (b) we get, $u_{p+1}=v_{p+1}=0$.
Hence the formal limit is the null solution.
Now, let us consider (\ref{1}), then, by using the previous notations and by
denoting \\
\begin{equation}
M(x,D)=OP \left\{-
\begin {pmatrix}
 0 & M^{12}(x,\xi)\\ M^{21}(x,\xi) & M^{22}(x,\xi)\\
\end{pmatrix} \right \},
\label{3}
\end{equation}
\\ we can rewrite the system (\ref{2}) in the following form\\
\begin{equation}
Z_t-M(x,D)Z=Q(x,Z^I,Z^{II})+D(Z^I,Z^{II}),
\label{4}
\end{equation}
where $ D(Z^I,Z^{II})=(D^I(Z^I,Z^{II}),D^{II}(Z^I,Z^{II}))$.\\  
We formulate  here  
the hypothesis concerning the existence of a symmetrizer for the
system (\ref{4}) in pseudodifferential form (see Taylor \cite{Tay81}), namely
\begin{itemize}
\item[\bf(A.4)] there exists $R(x,D)\! \in\! OPS^0 _{1,0}$ such that
$R(x,D)M(x,D)\!+(\!R(x,D)M(x,D)\!)^\ast\!\! \in\!\! OPS^0 _{1,0} $ and its symbol $R(x,\xi)$
is a positive definite matrix for $|\xi|>1$.
\end{itemize}
The next structure condition regards the existence of a symmetrizer for
the system (\ref{4}). We assume here a special ``block structure''
which is natural for 
strictly hyperbolic systems (for instance, in a more complicated 
framework, see  the seminal paper of 
Kreiss \cite{K70}, Majda and Osher 
\cite{MajOs} or Ralston \cite{Ral71}). The  block
structure follows also  for non - strictly hyperbolic systems having 
constant multiplicity, by a general result due to M\'etivier \cite{Met99}.
Here we are not assuming anyone of the previous conditions but directly the 
``block structure'' of the symmetrizer  $R(x,D)$.
\begin{itemize}
\item[\bf(S.2)]the symbol of $R(x,D)$ has the following form \\
$$R(x,\xi)= \begin{pmatrix} R_{11}(x,\xi) & 0\\ 0 & R_{22}(x,\xi)
\end{pmatrix}$$ \\where $ R_{11}(x,\xi)\in \mathcal{M}_{k\times k}$,
$R_{22}(x,\xi) \in \mathcal{M}_{(N-k)\times (N-k)} $ are symmetric
positive definite matrices and $R(x,D)M(x,D)+(R(x,D)M(x,D))^\ast \in
OPS^0 _{1,0} $.
\end{itemize}
Let us remark that in many applications this requirement will be 
automatically satisfyied.
\end{subsection}

\begin{subsection}{Dissipativity condition}
In this section we state the assumption on the nonhomgeneous term
$Q(x,Z^I,Z^{II})$.
\begin{itemize}
\item[\bf(D)]
$Q(x,Z^I,Z^{II})$ has the following form
$$Q(x,Z^I,Z^{II})=Q_0(x,Z^I,Z^{II})+Q_1(x,Z^I,Z^{II})$$ and
$Q(x,Z^I,0)=0$  for any $(x,Z^I)\in \ERRE ^d \times \ERRE ^k $
moreover \\ \\ 
{\bf(d1)} $Q_0(x,Z^I,Z^{II})\in C^1(\ERRE ^{N+d};
\ERRE ^{N-k})$, $Q_{0\nu}(x,Z^I,Z^{II})$ is bounded in \\$(x,Z^I,
Z^{II})$ and $[Q_{0 \nu}, R_{22}^{1/2}(x,D)]=0$. There exists
$\lambda _0>0$ such that for any $x\in \ERRE ^d$, $(Z^I, Z^{II})
\in \ERRE ^k \times \ERRE^ {N-k}$, $Q_{0 \nu}(x,Z^I, Z^{II})
\leq -\lambda _0I$,
\\ \\ {\bf(d2)} $Q_1(x,Z^I,Z^{II})\in
C^1(\ERRE ^{N+1}; \ERRE ^{N-k})$, $Q_{1\nu}(x,Z^I,Z^{II})$ is
bounded in \\$(x,Z^I, Z^{II})$ and the operator
$R_{22}(x,D)Q_{1\nu}$  satisfies \\
$\|R_{22}(x,D)Q_{1\nu}\|_{\mathcal{L}(L^2)}\leq \lambda _1$,
$\lambda _1>0$, $\lambda _1 \leq \lambda _0 /2$.
\end{itemize}
\begin{remark}
We splitted  the nonhomogeneous term $Q(x,Z^I,Z^{II})$ into
two terms that take into account  different dissipativity
mechanisms of $Q$. The former term $Q_0(\!x,\!Z^I,\!Z^{II}\!)$ cares about the
dissipativity of $Q$, and it commutes with $R_{22}(\!x,D\!)Q_{1\nu}$.
The latter  $Q_1(x,Z^I,Z^{II})$ does not commutes but defines
with $R_{22}(x,D\!)$ a bounded operator dominated by the
dissipative part of $Q_0(x,Z^I,Z^{II}\!)$.
\end{remark}
\begin{remark}
The class of dissipativity terms defined in (D) is not empty. In fact it
is sufficient to take $Q$ of the following form
\begin{equation}
Q(x,Z^I,Z^{II})=C(x)Z^{II} \label{33}
\end{equation}
with $C(x)\in {\mathcal M}_{(N-k)\times (N-k)}$ for any $x\in
\ERRE ^d$, $[R_{22}^{1/2}(x,D), C(x)]=0$  and there exists
$\gamma>0$ such that for any $x\in \ERRE ^d$, $C(x) \leq -\gamma
I$.
\end{remark}
\end{subsection}

\begin{subsection}{Formal analysis of the singular limit}

We will analyse the relaxation process of the following system

\begin{equation}
W_{t} +\frac{1}{\varepsilon} A(x,D)W = \frac{1}{\varepsilon ^2}B(x,
W)+\frac{1}{\varepsilon}D(W). \label{5}
\end{equation}
Following the construction of the previous paragraph we rewrite
 (\ref{5}) in this way\\
\begin{equation}
\left\{ \begin{array}{ll}
  Z^{I}_{t}+\displaystyle{\frac{1}{\varepsilon} M^{12}(x,D)Z^{II} =
  \frac{1}{\varepsilon}D^I(Z^I, Z^{II})}\\
  Z^{II}_{t}+\displaystyle{\frac{1}{\varepsilon} M^{21}(x,D)Z^{I}+
  \frac{1}{\varepsilon}M^{22}(x,D)Z^{II} = \frac{1}{\varepsilon ^2}
  Q(x,Z^I,Z^{II})+ \frac{1}{\varepsilon}D^{II}(Z^I, Z^{II})}.
  \end{array} \right.
\label{6}
\end{equation}\\
By formal asymptotics we are leaded to define
\begin{align}
Z^I(x, t) &= U^I(x,t) , & Z^{II}(x,t)&= \varepsilon U^{II}(x,t). &
\label{30}
\end{align}
The previous scaling has an equivalent interpretation as a scaling of
the time variable setting $\partial _\tau =\varepsilon \partial _t$,
for more details see \cite{MR}.  In this way the system
(\ref{6}) transforms into\\
\begin{equation}
\left\{ \begin{array}{ll} U^{I}_{t}+\displaystyle{M^{12}(x,D)U^{II} =
  \frac{1}{\varepsilon}D^I(U^I,\varepsilon U^{II})}\\ \\ \varepsilon
  ^2U^{II}_{t}+ \displaystyle {M^{21}(x,D)U^{I}+{\varepsilon}
  M^{22}(x,D) U^{II} = \frac{1}{\varepsilon} Q(x,U^I, \varepsilon
  U^{II})+D^{II}(U^I, \varepsilon U^{II})}.  \end{array} \right.
\label{7}
\end{equation}
If we denote by $(U^{I0}, U^{II0})$ the limit profile as $\varepsilon
\downarrow 0$, the formal limit system (\ref{7}) relaxes to
the system
\begin{equation}
\left\{ \begin{array}{ll} U^{I0}_{t}+\displaystyle{M^{12}(x,D)U^{II0}
  =D^I _ \nu (U^{I0},0) U^{II0}}\\ \\ \displaystyle {
  M^{21}(x,D)U^{I0} = Q_ \nu(x,U^{I0},0)U^{II0}+D^{II}(U^{I0}, 0)},
  \end{array} \right.
\label{8}
\end{equation}
where $Q_ \nu ,\ D^I _\nu, $ denote the derivative respect to the
variable $Z^{II}$respectively of $Q$ and $D^I$. By using the hypothesis 
(D), the
system (\ref{8}) transforms into $$
\left\{ \begin{array}{ll} U^{I0}_{t}+M^{12}(x,D) U^{II0} =D ^I _ \nu
(U^{I0},0) U^{II0}\\ \\U^{II0}= Q_ \nu^{-1}(x,U^{I0},0)
\left[M^{21}(x,D)U^{I0} -D^{II}(U^{I0}, 0)\right] \end{array} \right.
$$
That is equivalent to (by setting $U=U^{I0}$) the second order parabolic 
system\\
\begin{equation}\label{9}
 \begin{split}
U_{t}&+M^{12}(x,D)Q_ \nu^{-1}(x,U,0) M^{21}(x,D) U = M^{12}(x,D)Q_
\nu^{-1}(x,U,0) D^{II}(U,0) \\ & +D^I_ \nu(U,0)Q_
\nu^{-1}(x,U,0) \left[M^{21}(x,D)U - D^{II}(U, 0)\right].
\end{split}
\end{equation}
In the next sextion we will find sufficient conditions in order to
justify rigorously this formal analysis.
\end{subsection}
\end{section}

\begin{section}{A priori estimates and convergence analysis}
In this section we consider our rescaled system\\
\begin{equation}
\left\{ \begin{array}{ll} U^{I}_{t}+\displaystyle{M^{12}(x,D) U^{II} =
  \frac{1}{\varepsilon}D^I(U^I,\varepsilon U^{II})}\\ \varepsilon
  ^2U^{II}_{t}+ \displaystyle {M^{21}(x,D)U^{I}+{\varepsilon}
  M^{22}(x,D) U^{II} = \frac{1}{\varepsilon} Q(x,U^I, \varepsilon
  U^{II})+D^{II}(U^I, \varepsilon U^{II})}.  \end{array} \right.
\label{10}
\end{equation}
and we are going to develop the rigorous theory in order to get the
relaxed system (\ref{8}). We want to show that, as $\varepsilon
\downarrow 0$, the solutions of the rescaled system satisfy
\begin{align*}
U^{I} & \longrightarrow U^{I0} & &\text{strongly in
$L^{2}_{\loc}(\ERRE ^d\times \ERRE _+)$,} \\ U^{II} & \rightharpoonup
U^{II0}& &\text{weakly in $L^{2}(\ERRE ^d \times [0,T])$,} \\
\varepsilon U^{II} & \longrightarrow 0 & &\text{ strongly in
$L^{2}_{\loc}(\ERRE ^d \times\ERRE _+)$. }
\end{align*}
The basic idea used in this section is the assumption of the existence
of a symmetrizer for the system (\ref{2}) whose symbol has a suitable
block structure (conditions (A.4) and (S.2)).  The pseudodifferential
nature of the symmetrizer, as we will see later, will allow us to use
the properties of such kind of operators (Theorem 2.1, Theorem 2.2) in
order to obtain ``energy'' type estimate.
\begin{subsection}{A priori estimates}
In this section we wish to establish a priori estimates, independent of
$\varepsilon$, for the solution of the system (\ref{10}). To achieve this 
goal the following hypotheses are needed
\begin{itemize}
\item[\bf(B.1)] $M^{12}(x,\xi),M^{21}(x,\xi),M^{22}(x,\xi) \in S^1$,
\item[\bf(B.2)] $det \left[\left( M^{21}(x,\xi)\right)^T
M^{21}(x,\xi)\right]\neq 0$,
\item[\bf(B.3)]$D=(D^I(Z^I, Z^{II}), D^{II}(Z^I,Z^{II})) \in C^1(\ERRE
^N;\ERRE ^{k} \times \ERRE ^{N-k}) $, $D^{I}(Z^I,0)=0$,
$D^{II}(0,0)=0$, and $D^{I}_{\nu}$ is bounded in $(Z^I, Z^{II})$,
$D^{II}$ is a lipschitz function in $(Z^{I}, Z^{II})$, 
with lipschitz norm $\alpha$.
\end{itemize}
\begin{remark}
Since $M^{21}(x,\xi)\in \mathcal{M}_{(N-k)\times
k}$ from elementary linear algebra we deduce condition (B.2) is
violated whenever $k>\frac{N}{2}$.
\end{remark}
The next results concerns the ``energy'' type estimates independent on 
$\varepsilon$, obtained via the existence of symmetrizers and via
the dissipativity conditions (D).
\begin{theorem}
Let us consider the solution $\{U^I\}$, $\{U^{II}\}$ of the Cauchy
problem for the scaled system (\ref{10}). 
Assume $U ^I(x,0) \in \left[L^2(\ERRE ^d) \right]^k $,
 $U^{II}(x,0) \in \left[L^2(\ERRE ^d) \right]^{N-k}$ and 
the hypotheses (A.4), (S.2), (B.1),
(B.2), (B.3), (D) hold. Then, there exists $\varepsilon_{0}>0$, such that 
for $\varepsilon \in (0,\varepsilon_{0})$,
one has
\begin{itemize}
\item[\bf(i)] for all $T>0$,  there exists
$M(T)>0$, independent from $\varepsilon$, such
 that\\
$\|U^{II} \|_{L^2(\ERRE ^d \times [0,T])} \leq M(T)$ and
 $\displaystyle \sup _{[0,T]}
\|\varepsilon U^{II} (\cdot,t)\|\leq M(T)$,
\item[\bf(ii)]$\{\varepsilon^2U^{II} _t\}\quad \text{ is relatively compact in
 $H^{-1}_{loc}(\ERRE ^d \times \ERRE_+)$} $,
\item[\bf(iii)] $\{U^I\}$ is uniformely bounded, with respect to
$\varepsilon$, in $ L^\infty \left(\ERRE_+,L^{2}(\ERRE ^d )\right)
$, namely for all $T>0$,  there exists $M(T)>0$, independent from
$\varepsilon$, such that \\$\displaystyle \sup _{[0,T]} \|U^I
(\cdot,t)\|\leq M(T)$.
\end{itemize}
\label{t1}
\end{theorem}
\begin{proof}
To simplify our notation, we set
\begin{align*}
 U &=\begin{bmatrix}
    U^I \\
    U^{II}
\end{bmatrix}, &
U^\varepsilon &=\begin{bmatrix}
    U^I \\
    \varepsilon U^{II}
\end{bmatrix}, &
F &=\begin{bmatrix}
    \displaystyle{\frac{1}{\varepsilon}D^I(U^I,\varepsilon U^{II})} \\
    \displaystyle{\frac{1}{\varepsilon}
Q(x,U^I, \varepsilon U^{II})+D^{II}(U^I, \varepsilon U^{II})}
\end{bmatrix}
\end{align*}
$$ M^\varepsilon (x,D)=OP \left\{-
\begin {pmatrix}
 0 &  M^{12}(x,\xi)\\
M^{21}(x,\xi ) &  \varepsilon M^{22}(x,\xi )\\
\end{pmatrix} \right \}
$$ 
By using  the Theorem \ref{p2}, let  $R(x,D)$ be the symmetrizer given
by the hypotheses  (A.4), (S.2), then there exists $\widetilde{R}=R \ mod\
OPS^{-\infty}$, with $\widetilde{R}\geq \eta I>0$. So we can write
$\widetilde{R}=R+T$, $T\in OPS^{-\infty}$, where by the block
structure of $R(x,D)$, $T(x,\xi)$ takes the following form
 $T(x,\xi)=\left(\begin{smallmatrix}
    T_{11}(x,\xi) & 0   \\
    0 & T_{22}(x,\xi)
    \end{smallmatrix}\right)$. 
Hence it follows    
\begin{align*}
\frac{d}{dt}(\widetilde{R}(x,D)U^\varepsilon, U^\varepsilon)_2&\!=
(\widetilde{R}(x,D)U^\varepsilon _t, U^\varepsilon)_2 +
(\widetilde{R}(x,D)U^\varepsilon, U^\varepsilon _t)_2\\
&\!=\!(\!\widetilde{R}(x,D)(M^\varepsilon(x,D)U+F\!),
U\!)_2\!+\!(\!\widetilde{R}(x,D)U,M^\varepsilon(x,D)U \!+\!F)_2\\
&\!=\!((\widetilde{R}(x,D)M^\varepsilon(x,D)+(\widetilde{R}(x,D)M^\varepsilon(x,D))^\ast)
U , U)_2\\&\quad +2 (\widetilde{R}(x,D)F,U)_2 =I_1 +I_2 .
\end{align*}
We estimate separately $I_1$ and $I_2$. By applying Theorem \ref{p1}
and Theorem \ref{p2} and by using the hypothesis (A.4) for $I_1$ we get
$$ I_1 \leq 2 \|U^{II}\|\|U^I\| +\varepsilon \|U^{II}\|^2.$$ Let us
focus our attention on $I_2$, then one has
\begin{align*}
I_2 &=(\widetilde{R}(x,D)D(U^I, \varepsilon U^{II}),U)_2+
(\widetilde{R}(x,D)\begin{pmatrix}
            0 \\
    \varepsilon ^{-1}
Q(x,U^I, \varepsilon U^{II})
\end{pmatrix} ,U)_2
\\
&=(\varepsilon ^{-1}D^I(U^I, \varepsilon
U^{II}),(R_{11}(x,D)+T_{11}(x,D))U^I)_2\\ &\quad+ (D^{II}(U^I,
\varepsilon U^{II}),(R_{22}(x,D)+T_{22}(x,D))U^{II})_2\\ &\quad +
(\varepsilon ^{-1}R_{22}(x,D)Q(x,U^I, \varepsilon
U^{II}),U^{II})_2+ (\varepsilon ^{-1}T_{22}(x,D) Q(x,U^I,
\varepsilon U^{II}) ,U^{II})_2\\ &=I_{2,1}+I_{2,2}+I_{2,3}+I_{2,4}.
\end{align*}
To estimate $I_{2,1}$ and $I_{2,2}$ we use hypothesis (B.4) and we get
$$I_{2,1}+I_{2,2}\leq \alpha\|U^{II}\|\|U^I\|
+\alpha\varepsilon \|U^{II}\|^2 .$$ Now we turn to $I_{2,3}$.
Let us denote by $\bar Q _ \nu = \int _0 ^1 Q_\nu(x,U^I, \varepsilon
\theta U^{II})d\theta $, then  we can rewrite $I_{2,3}$ in the following
way
\begin{align*}
I_{2,3}&=(R_{22}(x,D)\bar Q _ {0\nu} U^{II},U^{II})_2 + (
R_{22}(x,D)\bar Q _{1\nu}  U^{II}, U^{II})_2=I_{2,31}+I_{2,32}.
\end{align*}
From (d1), it follows 
\begin{align*}
I_{2,31}&=\!(R_{22}^{1/2}(x,D)\bar Q _{0\nu}
U^{II}, R_{22}^{1/2}(x,D)U^{II})_2\!=\!(\!\bar Q
_{0\nu}R_{22}^{1/2}(x,D) U^{II}, R_{22}^{1/2}(x,D)U^{II}\!)_2\\&\leq
-\lambda _0 \|U^{II}\|^2.
\end{align*}
By using (d2) we get
\begin{align*}
I_{2,32}&=(R_{22}^{1/2}(x,D)\bar Q _{1\nu} U^{II},
R_{22}^{1/2}(x,D)U^{II})_2  \leq \lambda _1 \|U^{II}\|^2 \leq
\frac{\lambda _0}{2} \|U^{II}\|^2,
\end{align*}
hence $$I_{2,3}=- \frac{\lambda _0}{2} \|U^{II}\|^2 .$$ Now
it remains to estimate $I_{2,4}$. For any $\delta _1>0$ we have
$$I_{2,4}\leq \delta _1 \|U^{II}\|^{2}+\frac{c'}{\delta _1}
\|U^{II}\|^{2}_{H^{-1}}$$ 
By adding $I_1$, $I_2$, for all 
$\delta _1 >0,\ \delta _2 >0$, we get
\begin{align*}
\frac{d}{dt}(R(x,D)U^\varepsilon, U^\varepsilon)_2
&\leq (2+\alpha)\|U^{II}\|\|U^I\| +\varepsilon(\alpha+1)
\|U^{II}\|^2 -\frac{\lambda _0}{2}\|U^{II}\|^2\\&\quad+\delta _1 
\|U^{II}\|^{2}+\frac{c'}{\delta _1}
\|U^{II}\|^{2}_{H^{-1}}
 \\
&\leq (\delta _1 +\delta _2) \|U^{II}\|^2+\frac{2+\alpha}{\delta
_2 }\|U^{I}\|^2 +\varepsilon(\alpha+1) \|U^{II}\|^2 \\&\quad-\frac{\lambda
_0}{2}\|U^{II}\|^2 +\frac{c'}{\delta _1} \|U^{II}\|^{2}_{H^{-1}}.
\end{align*}
If we  choose $\displaystyle{\varepsilon
<\frac{\lambda _0 }{8(\alpha+1)}}$ and $\displaystyle{(\delta _1
+\delta _2) <\frac{\lambda _0}{8}}$ it follows 
\begin{equation}
\frac{d}{dt}(R(x,D)U^\varepsilon, U^\varepsilon)_2 \leq
-\frac{\lambda _0}{4}\|U^{II}\|^2+c\|U^{I}\|^2 +
c\|U^{II}\|^{2}_{H^{-1}}. \label{11}
\end{equation}
Let us denote by
\begin{equation}
E(t)=\int_{\ERRE ^d}\varepsilon^2|U^{II}(x,t)|^2dx+
\int_{\ERRE ^d} \left|U^I(x,t)\right|^2dx.
\label{19}
\end{equation}
By the second equation of the system (\ref{10}) and by the
smoothness of the coefficients we get the following estimate
$$\|U^{II}\|^2_{H^{-1}}\leq c \sqrt {E(t)}.$$ By integrating
(\ref{11}) on $[0,t] $ we obtain the energy $E(t)$ satisfies
the following inequality
\begin{equation}
E(t)\leq E(0)-\frac{\lambda _0}{4}\int_{0}^{t}\!\!\int_{\ERRE ^d
}|U^{II}(x,t)|^2dxds +c\int_{0}^{t}E(s)ds \label{20}
\end{equation}
for some  $c>0$, then by applying the 
Gronwall's lemma, we have
\begin{equation}
E(t)\leq E(0) e^{ct},
\end{equation}
and 
\begin{equation}
\int_{0}^{t}\!\!\int_{\ERRE ^d}|U^{II}|^2dxds \leq
E(0)\left(e^{ct}+2\right).
\end{equation}
We can  conclude that for any $T>0$ there exists $M(T)>0$, indipendent
from $\varepsilon$, such that
\begin{align*}
\|U^{II}\|_{L^2(\ERRE ^d \times[0,T])} & \leq M(T), & \sup _{[0,T]}\|U^{II} (\cdot,t)\|& \leq M(T),&  \sup _{[0,T]}\|U^{I}
(\cdot,t)\|& \leq M(T).
\end{align*}
In this way we proved (i) and (iii). Let us consider $\omega$
relatively compact subset of $\ERRE ^d \times \ERRE_+$, then
\begin{align*}
\|\varepsilon^2U^{II}  _t\|_{H^{-1}(\omega)}&=\sup
_{\|\phi\|_{H^{1}_{0}(\omega)}=1}|\langle \varepsilon^2 U^{II} _t
, \phi\rangle| =\sup
_{\|\phi\|_{H^{1}_{0}(\omega)}=1}\left|\int\!\!\int \varepsilon^2
U^{II}\phi _t dxdt\right|\\ & \leq  \varepsilon^2 \sup
_{\|\phi\|_{H^{1}_{0}(\omega)}=1}(\|U^{II}\|  \|\phi _t\|) \leq
M(T) \varepsilon ^2.
\end{align*}
\end{proof}
\begin{remark}
If there exists a way to control $\| Q _ \nu\|_{Lip}$ we do not need any
assumption (D) but to
estimate $I_{2,3}$ it is sufficient to have for any
$(x,Z^I, Z^{II}) \in \ERRE ^d \times \ERRE ^k \times \ERRE^
{N-k}$, $Q_ \nu(x,Z^I, Z^{II}) \leq -\lambda I$, $\lambda >0$.
Indeed 
\begin{align*}
I_{2,3}&=(R_{22}^{1/2}(x,D)\bar Q _ \nu U^{II}, R_{22}^{1/2}(x,D)U^{II})_2\\
&=(\bar Q _ \nu R_{22}^{1/2}(x,D) U^{II},R_{22}^{1/2}(x,D) U^{II})_2
        +([R_{22}^{1/2}(x,D),\bar Q _ \nu]U^{II}, U^{II})_2,
\end{align*}
then by taking into account the Proposition \ref{p3} and the dissipative
condition on $ Q_\nu(x,U^I, \varepsilon U^{II})$, we have
$$I_{2,3}\leq -\frac{\lambda}{2}\|U^{II}\|^2
+c\|U^{II}\|^2_{H^{-1}}. $$ which lead to the energy inequality
(\ref{20}).
\end{remark}
\end{subsection}

\begin{subsection}{Strong convergence analysis}
We begin with an immediate
consequence of (i) and (ii) in the  Theorem (\ref{t1}).
\begin{theorem}
Let \ us \ consider \ the \ solution \ $\{U^{II}\}$ \ of \ the \ Cauchy \
problem\  for\  system (\ref{10}). Assume\    
$U ^I(x,0) \in \left[L^2(\ERRE ^d) \right]^k $,
 $U^{II}(x,0) \in \left[L^2(\ERRE ^d) \right]^{N-k}$ and moreover  
the hypotheses  (A.4), (S.2), (B.1), (B.2),
(B.3), (D) hold. Then there exists $U^{II0}\in
[L^2(\ERRE ^d\times [0,T])]^{N-k}$, such that, as $\varepsilon
\downarrow0$, one has (extracting eventually a subsequence)
\begin{align}
& U^{II} \rightharpoonup U^{II0} \hspace{3mm}& &\text{ weakly
in $L^{2}(\ERRE ^d\times [0,T])$} \label{12}\\
&\varepsilon
U^{II}\longrightarrow 0 \hspace{-4mm}& & \text{ strongly in
$L^{2}_{loc}(\ERRE ^d \times \ERRE_+)$} \label{13}\\
&\{\varepsilon^2
U^{II} _t\}\longrightarrow 0 \hspace{-4mm}& & \text{  in
$H^{-1}_{loc}(\ERRE ^d \times \ERRE_+)$}. \label{14}
\end{align}
\label{t2}
\end{theorem}
Our next step is to prove the strong  convergence for the sequence
$\{U^{I}\}$ in the norm of $L^{2}_{loc}(\ERRE ^d \times \ERRE_+)$.  For 
this porpuse
we only  use the estimates obtained in the previous
section. Our main tool in the limit process will be  
Tartar's  and G\'erard's Theorem \ref{g}, \cite{TarM}, \cite{Ger}.
\begin{theorem}
Let\  us\  consider\  the\  solution \ $\{U^I\}$,
 of\ the\  Cauchy\  problem\  for\  system (\ref{10}).
Assume\   $U ^I(x,0) \in \left[L^2(\ERRE ^d) \right]^k $,
 $U^{II}(x,0) \in \left[L^2(\ERRE ^d) \right]^{N-k}$ and the hypotheses  (A.4), (S.2), (B.1), (B.2), (B.3),(D) hold.
Then there exists $U^{I0}\in[ L^2(\ERRE ^d\times [0,T])]^k$, such that,
as $\varepsilon \downarrow0$, one has (extracting eventually
subsequences)
\begin{align}
& U^I \longrightarrow U^{I0} &  &\text{ strongly in
$L^{2}_{loc}(\ERRE ^d \times \ERRE_+)$}  \label{15}
\end{align}
\label{t3}
\end{theorem}
\begin{proof}
By using the hypothesis (D), $\bar Q _ \nu(x,U^I,\!\varepsilon
U^{II})$ \ is uniformely bounded  in $L^\infty$ then
$\varepsilon ^{-1}Q _ \nu(x,U^I,\!\varepsilon U^{II})= \bar Q _
\nu(x,U^I,\varepsilon U^{II})U^{II}$ is uniformely bounded in
$L^2$, therefore $\varepsilon ^{-1}Q _ \nu(x,U^I,\varepsilon
U^{II})$ is relatively compact in $H_{loc}^{-1}$. In a similar
way, thanks to the conditions on the function $D(U^I,\varepsilon
U^{II})$, we get $\varepsilon ^{-1}D^I(U^I,\varepsilon U^{II})$
and $D^{II} (U^I,\varepsilon U^{II})$  relatively compact in
$H_{loc}^{-1}$. Now by using the identities $$
\begin{array}{ll}
 U^{I}_{t}+M^{12}(x, D)U^{II} =\displaystyle{
\frac{1}{\varepsilon}D^I(U^I,\varepsilon U^{II})}\\
\displaystyle
{M^{21}(x, D)U^{I}= \frac{1}{\varepsilon} Q(x,U^I, \varepsilon
U^{II})+D^{II}(U^I, \varepsilon U^{II})}- \varepsilon
^2U^{II}_{t}-{\varepsilon} M^{22}(x, D)U^{II}
\end{array}
$$
and by taking into account the a priori estimates and the smoothness of the
coefficients of the system, we can conclude
\begin{equation}
    \begin{pmatrix}
     U^{I}_{t}+\displaystyle{M^{12}(x,D)U^{II}} \\
    \displaystyle {M^{21}(x,D)U^{I}}
    \end{pmatrix} \quad \quad \text{is relatively compact in $(H^{-1}_{loc})^2$.}
    \label{100}
\end{equation}
Since the symbols $M^{ij}(x,\xi)$ are polihomogeneous in the sense of the 
Section 2, we can decompose the operator $M^{ij}(x,D)$ in the following way,
\begin{equation}
    M^{ij}(x,D)=M^{ij}_{1}(x,D)+M^{ij}_{r}(x,D),
    \label{101}
\end{equation}
where $M^{ij}_{1}(x,D)\in OPS^{1}$, while $M^{ij}_{r}(x,D)\in OPS^{0}$. By 
applying the Theorem \ref{p1} and by the decomposition (\ref{101})
we get
\begin{equation}
    \begin{pmatrix}
     U^{I}_{t}+\displaystyle{M^{12}_{1}(x,D)U^{II}} \\
    \displaystyle {M^{21}_{1}(x,D)U^{I}}
    \end{pmatrix} \quad \quad \text{is relatively compact in $(H^{-1}_{loc})^2$.}
    \label{102}
\end{equation}
In order to fit into the framework of the Theorem \ref{g} we set
$$P\begin{bmatrix}
U^I \\
U^{II}
\end{bmatrix}= \begin{bmatrix}
                 I_{k\times k} & 0\\
                          0    & 0\\
\end{bmatrix}
\partial _t\begin{bmatrix} U^I \\
                         U^{II}\\
\end{bmatrix} +\displaystyle{ \begin{bmatrix} 0 & M^{12}_{1}(x,D) \\
               M^{21}_{1}(x,D) & 0
\end{bmatrix}}\begin{bmatrix}
                          U^I \\
                         U^{II}\\
\end{bmatrix}.$$
The principal symbol of $P$ is given by
$$p(x,\xi)=\begin{bmatrix}
                I_{k\times k} &  0\\
                          0   & 0\\
\end{bmatrix} \xi_0 +
\displaystyle{\begin{bmatrix} 0 & M^{12}_{1} (x, \xi') \\
               M^{21}_{1}(x, \xi ') & 0 \\
\end{bmatrix}},$$ for any $\xi=(\xi ',\xi _0)\in \ERRE^{d+1}$, $|\xi|=1$.
We
notice that
$$p(x,\xi)\begin{bmatrix}\lambda \\
                 \mu\\
\end{bmatrix}=0 \Longleftrightarrow \left\{\begin{array}{ll}
       \xi_0\lambda+
M^{12}_{1}(x, \xi ') \mu=0\\  M^{21}_1 (x,\xi ') \lambda=0
       \end{array}
\right. $$ for all $\lambda\in \ERRE^k$, $\mu \in \ERRE ^{N-k}$.\\
Now if $\xi '=0$ then $\xi_0 \neq 0$ and so $\lambda =0$,
otherwise if $\xi ' \neq 0$ by using the hypothesis
(B.2) we get $\lambda =0$.Therefore we get 
$$\left\{ (x,\xi,\lambda,\mu)\quad
\text{such that}\quad p(x,\xi)\begin{bmatrix}\lambda \\
                 \mu\\
\end{bmatrix}=0 \right\}\subset \left\{ \lambda ~|~ \lambda =0\right\}.$$
We take now
$$q(x)=\begin{bmatrix}I_{k \times k} & 0\\
           0 & 0 \\
           \end{bmatrix}$$
and for all $\xi \neq 0$, $\xi=(\xi_0, \xi ')$ we have
$$p(x,\xi)\begin{bmatrix}\lambda \\
                 \mu\\
\end{bmatrix}=0 \qquad
\text{implies} \qquad \langle q(x)\begin{bmatrix}\lambda\\
                 \mu\\
\end{bmatrix}, \begin{bmatrix}\lambda \\
                 \mu\\
\end{bmatrix}\rangle =0$$ for all $\lambda\in \ERRE^k$, $\mu \in \ERRE
^{N-k}$.\\ Now \ we \  can \  apply  \ Theorem \ (\ref{g}) of G\'erard and we
conclude that for any $\varphi\in \mathcal{D}(\ERRE ^d\times \ERRE_+)$ \\ 
$$\int
_{0}^{+ \infty}\!\!\int _{\ERRE ^d}
\left|U^I(x,t)\right|^2\varphi(x,t)dxdt \longrightarrow
\int_{0}^{+ \infty}\!\!\int _{\ERRE ^d}
\left|U^{I0}(x,t)\right|^2\varphi(x,t)dxdt$$ where $U^{I0}$
denotes, in view of Theorem (\ref{t1}) the weak limit of $U^I$ in
$L^2(\ERRE ^d\times \ERRE_+)$. In this way we obtain\\ $$ U^I
\longrightarrow U^{I0} \qquad \text{strongly in $L^{2}_{loc}(\ERRE
^d \times \ERRE_+)$}. $$
\end{proof}
\begin{corollary}
Assume that the hypotheses of Theorems (\ref{t2}), (\ref{t3})
 hold, then $(U^{I0},U^{II0})$ verifies the following system,
 in the sense of distributions, \\
\begin{equation}
\left\{
  \begin{array}{ll}
    U^{I0}_{t}+M^{12}(x,D)U^{II0}
=D^I _ \nu (U^{I0},0) U^{II0}\\ \\ M^{21}(x, D)U^{I0} =Q_
\nu(x,U^{I0},0)U^{II0}+D^{II}(U^{I0}, 0).
  \end{array}
\right.
\label{16}
\end{equation}
\end{corollary}
\begin{proof}
By
taking into account the regularity hypothesis on $Q$, the strong convergence
of $\{ U^I,\varepsilon U^{II} \}$ in $L^2 _{loc}$ and the weak convergence of
$U^{II}$, we get \\
$$ \varepsilon ^{-1}Q _ \nu(x,U^I,\varepsilon U^{II}) \rightharpoonup
 Q  _\nu(x,U^{I0},0)U^{II0} \qquad \text {weakly in
$L^{2}_{loc}(\ERRE ^d \times \ERRE_+)$}.$$
In a similar way we have \\
$$ \varepsilon ^{-1}D^I(U^I,\varepsilon U^{II}) \rightharpoonup
 D^I_ \nu(U^{I0},0)U^{II0} \qquad \text {weakly in
$L^{2}_{loc}(\ERRE ^d \times \ERRE_+)$}.$$
Therefore  we can pass into the limit the other terms and we obtain the 
relations (\ref{16}).
\end{proof}
\begin{remark}
We can weakly relax the assumption concerning the polihomogeneity of our 
symbols by assuming directly the decomposition (\ref{101})
\end{remark}
\end{subsection}

\begin{subsection}{Parabolicity condition}

Let us restrict our attention to the differential operator case, namely 
(\ref{0}) takes the form
$$W_{t} +\frac{1}{\varepsilon} \sum_{j =
1}^{d}A_{j}(x)\partial_{j}W = \frac{1}{\varepsilon ^2}B(x,
W)+\frac{1}{\varepsilon}D(W).$$
In this case the identities (\ref{16}) become
\begin{equation}
\left\{
  \begin{array}{ll}
    U^{I0}_{t}+\displaystyle{\sum_{j=1}^{d}M^{12}_{j}(x)\partial _j U^{II0}
=D^I _ \nu (U^{I0},0) U^{II0}}\\
 \displaystyle {   \sum_{j=1}^{d}M^{21}_{j}(x)\partial _j U^{I0} =
   Q_ \nu(x,U^{I0},0)U^{II0}+D^{II}(U^{I0}, 0)},
  \end{array}
\right. \label{161}
\end{equation}
which is equivalent to write (where we set
$U=U^{I0}$)\\
\begin{equation}
\begin{split}
 U_{t}&+\sum_{j=1}^{d}M^{12}_{j}(x)\partial _j\left(Q_ \nu^{-1}(x,U,0)
\sum_{k=1}^{d} M^{21}_{k}(x)\partial _{k} U \right)\\&=
\sum_{j=1}^{d}M^{12}_{j}(x)\partial_j(Q_ \nu^{-1}(x,U,0) D^{II}(U,0)) \\
 &+D^I_ \nu(U,0)Q_ \nu^{-1}(x,U,0)
\left[\sum_{k=1}^{d}M^{21}_{k}(x)\partial _k U - D^{II}(U, 0)\right].
\end{split}
\label{26}
\end{equation}
We want to show, in this simpler case, that (\ref{26}) is parabolic
in the sense of Petrowski as recalled in the Section 2. 
Taking into account the notations
of Section 2,  one has $$L=M^{12}(x,\xi)Q_ \nu^{-1}(x,U,0) M^{21}(x,\xi)
\in S^2_{1,0},$$ where $M^{12}(x,\xi)=\displaystyle{i\sum _{j=1}^d
M^{12}_j(x) \xi _j}$ and $M^{21}(x,\xi)=\displaystyle{i\sum
_{j=1}^d M^{12}_j(x) \xi _j}$. We are going to show the
existence of a positive matrix $P_0$ such that $P_0 L+L^\ast P_0$
is negative definite. By the definition of symmetrizer for the
system (\ref{2}) we know the following relations between the blocks of
$R(x,\xi)$ and the coefficients of the system (\ref{2}) 
\begin{equation} 
\label{17}
\begin{split}
 R_{11}(x,\xi)M^{12}(x,\xi)&=(M^{21}(x,\xi))^TR_{22}(x,\xi), \\
  R_{22}(x,\xi)M^{21}(x,\xi)&=(M^{12}(x,\xi))^TR_{11}(x,\xi).
\end{split}
\end{equation}
Now let us denote by
$$P_0= R_{11}(x,\xi), $$
we have to prove
$$P_0 L +L^\ast P_0 \leq -\beta I \qquad \qquad \beta >0 .$$
By using the relations (\ref{17}) and the condition (D) it follows
\begin{align*}
\left((P_0 L +L^\ast P_0)\eta, \eta  \right)&=
\left(R_{22}(x,\xi)M^{12}(x,\xi)Q_
\nu^{-1}(x,U,0)M^{21}(x,\xi)\eta, \eta \right)\\&+
\left(M^{21}(x,\xi))^T (Q_
\nu^{-1}(x,U,0))^T(M^{12}(x,\xi))^TR_{11}(x,\xi)\eta
,\eta\right)\\ &=\left((M^{21}(x,\xi))^T(R_{22}(x,\xi)Q_
\nu^{-1}(x,U,0)\eta, \eta \right)\\ &+\left((Q_
\nu^{-1}(x,U,0))^TR_{22}(x,\xi))M^{21}(x,\xi)\eta, \eta \right)\\
&\leq -\lambda _0\left|M^{21}(x,\xi)\eta \right|^2 \leq -\beta
\left|\eta \right|^2 
\end{align*} with $\beta >0$, this is the notion of parabolicity introduced 
in the Section 2.
\end{subsection}
\end{section}

\begin{section}{The constant coefficients case}
In this section we want to show how in the case of constant coefficient
differential semilinear systems our theory can be easily simplified. 
First we remark that in this case we don't need to use
pseudodifferential theory because of the constant coefficient we
can handle them with classical methods. We point out that also in
this case we will assume the existence of  symmetrizers $R(x,D)$
with block structure but since the coefficients are constant the
principal symbol of $R(x,D)$ depends only on the variable $\xi$
more exactly it is a homogeneous radial function  of degree zero of
$\xi$, $R(\xi)=R\left(\frac{\xi}{|\xi|}\right)$. Then the symmetrizers 
reduce to Fourier multipliers.
\begin{subsection}{A priori estimate}
Here we consider the following system\\
\begin{equation}
\left\{
  \begin{array}{ll}
    U^{I}_{t}+\displaystyle{\sum_{j=1}^{d}M^{12}_{j}\partial _j U^{II} =
\frac{1}{\varepsilon}D^I(U^I,\varepsilon U^{II})}\\
    \varepsilon ^2U^{II}_{t}+
\displaystyle {\sum_{j=1}^{d}M^{21}_{j}\partial_j U^{I}+{\varepsilon}
\sum_{j=1}^{d}M^{22}_{j}\partial _j U^{II} =
\frac{1}{\varepsilon}
Q(x,U^I, \varepsilon U^{II})+D^{II}(U^I, \varepsilon U^{II})}.
  \end{array}
\right.
\label{18}
\end{equation}
with the same hypotheses  (A.4), (S.2), (B.1),  (B.2), (B.3),  (D) of
Section 3.2 and  Section 4.1, specialized to our simpler framework. 
Then we have also in this case the following theorem
\begin{theorem}
Let us consider the solution $\{U^I\}$, $\{U^{II}\}$ of the Cauchy
problem for the system (\ref{18}). Assume 
$U ^I(x,0) \in \left[L^2(\ERRE ^d) \right]^k $,
$U^{II}(x,0) \in \left[L^2(\ERRE ^d) \right]^{N-k}$ and that  
hypotheses (A.4), (S.2), (B.1),
(B.2), (B.3), (D) hold. Then there exists $\varepsilon_{0}>0$,
such that for $\varepsilon \in (0,\varepsilon_{0})$,
one has
\begin{itemize}
\item[{\bf(i)}] for any $T>0$,  there exists $M(T)>0$, independent
from $\varepsilon$, such that \\$\|U^{II} \|_{L^2(\ERRE ^d \times
[0,T])} \leq M(T)$ and
 $\displaystyle \sup _{[0,T]}
\|\varepsilon U^{II} (\cdot,t)\|\leq M(T)$,
\item[{\bf(ii)}]$\{\varepsilon^2U^{II} _t\}\quad 
\text{ is relatively compact in
 $H^{-1}_{loc}(\ERRE ^d \times \ERRE_+)$} $,
\item[{\bf(iii)}]$\{U^I\}$ is uniformely bounded, with respect to
$\varepsilon$, in $ L^\infty \left(\ERRE_+,L^{2}(\ERRE ^d )\right)
$, namely for any $T>0$, there exists $M(T)>0$, independent from
$\varepsilon$, such that\\ $\displaystyle \sup _{[0,T]} \|U^I
(\cdot,t)\|\leq M(T)$.
\end{itemize}
\label{t4}
\end{theorem}
\begin{proof}
We apply the Fourier transform to the system (\ref{18}) and we
multiply in $L^2(\ERRE ^d)$ by $R(\xi)\widehat {U}$ obtaining the following
inequality
\begin{align*}
\frac{d}{dt}\left\{ \frac{\varepsilon ^2}{2}
|R_{22}^{1/2}(\xi)\widehat{U^{II}}|^2 +
|R_{11}^{1/2}(\xi)\widehat{U^I}|^2
\right\}&=(R_{11}(\xi)\widehat{U^I}, \varepsilon
^{-1}\widehat{D^I}(U^I,\varepsilon U^{II})_2 \\ &\hspace{-4cm}+(
R_{22}(\xi)\widehat{U^{II}}, \varepsilon ^{-1}\widehat
{Q}(U^I,\varepsilon U^{II} ))_2 +( R_{22}(\xi)\widehat{U^{II}},
\widehat{D^{II}}(U^I, \varepsilon U^{II}))_2 .
\end{align*}
Defining the energy as in (\ref{19}), taking into account the hypotheses,
the properties of the symmetrizer and Plancharel theorem it yields the
standard energy inequality (\ref{20}). The remaining part of the proof
follows exactly the same arguments used in the previous 
section so it
is omitted.
\end{proof}
\end{subsection}
\begin{subsection}{Basic ideas on strong convergence}
The analysis of strong convergence in this case reduces to analyse the 
convergence of quadratic forms, hence it can be obtained via the classical
compensated compactness  result of Tartar (see \cite{Tar79}, \cite{Tar83}, \cite{Mur78}
see also \cite{Dac82}).
\begin{theorem}\label{cc}
{\bf (Tartar's Compensated compactness)}\\
Let us consider
\begin{enumerate}
\item
a bounded open set $\Omega\subset\mathbb{R}^{n}$;
\item
a sequence $\{l^{\nu}\}_{\nu = 1}^{\infty}$,
$l^{\nu}:~\Omega\subset \mathbb{R}^{n}
\longrightarrow\mathbb{R}^{m} $;
\item
a symmetric matrix
$\Theta:~\mathbb{R}^{m}\longrightarrow\mathbb{R}^{m} $;
\item
constants $a_{jk}^{i}\in \mathbb{R}$, $i = 1,\ldots,q$,
$j = 1,\ldots,m$, $k = 1,\ldots,n$.
\end{enumerate}
Let us define
\begin{align*}
f(\alpha) & = \left\langle \Theta \alpha,\alpha \right\rangle,
\quad\text{for all $\alpha \in \mathbb{R}^{m}$;} \\
\Lambda & = \left\{ \lambda \in \mathbb{R}^{m}:~
\exists \xi \in \mathbb{R}^{n}\setminus
\{0\}, ~\sum\limits _{j,k}a_{jk}^{i} \lambda_{j}
\xi_{k} = 0, i = 1,\ldots,q\right\}.
\end{align*}
Assume that
\begin{itemize}
\item[\bf(a)]
there exists $\widetilde{l}\in L^{2}_{m}\left(\Omega\right)$
such that $l^{\nu}\rightharpoonup \widetilde{l}$ in
$L^{2}_{m}\left(\Omega\right)$ as $\nu\uparrow\infty$;
\item[\bf(b)]
$\mathcal{A}^{i}l^{\nu} = \sum\limits _{j,k}a^{i}_{jk}\frac{\partial
l_{j}^{\nu}}{\partial x_{k}}$, $i=1,\ldots q$ are relatively compact in
$H^{-1}_{\loc}\left(\Omega\right)$;
\item[\bf(c)]
$f_{|\Lambda}\equiv 0$;
\item[\bf(d)] there exists
$\widetilde{f}\in \mathbb{R}$ such that $f(l)\rightharpoonup
\widetilde{f}$ in the sense of measures $\mathcal{M}(\Omega)$.
\end{itemize}
Then we have $\widetilde{f} = f(\widetilde{l})$.
\end{theorem}
Now we can state our convergence result
\begin{theorem}
Let  $\{U^I, U^{II}\}$ be the solution of the Cauchy
problem for system (\ref{18}). Let us assume the hypotheses  
(A.4), (S.2),(B.1), (B.2),
(B.3), (B.4), (B.6) hold, then there exists
$U^{I0}\in [L^2(\ERRE ^d\times [0,T])]^{k}$ and
$U^{II0}\in [L^2(\ERRE ^d\times [0,T])]^{N-k}$,
such that, as $\varepsilon
\downarrow0$, (extracting eventually subsequences)
\begin{align}
& U^{II} \rightharpoonup U^{II0} \hspace{3mm}& &\text{ weakly
in $L^{2}(\ERRE ^d\times [0,T])$} \label{21}\\
&\varepsilon
U^{II}\longrightarrow 0 \hspace{-4mm}& & \text{ strongly in
$L^{2}_{loc}(\ERRE ^d \times \ERRE_+)$} \label{22}\\
&
U^{I} \longrightarrow  U^{I0} \hspace{-4mm}& & \text{ strongly in
$L^{2}_{loc}(\ERRE ^d \times \ERRE_+)$} \label{23}
\end{align}
and the limit profile $(U^{I0}, U^{II0})$ verifies the following system
\begin{equation}
\left\{
  \begin{array}{ll}
    U^{I0}_{t}+\displaystyle{\sum_{j=1}^{d}M^{12}_{j}\partial _j U^{II0}
=D^I _ \nu (U^{I0},0) U^{II0}}\\
 \displaystyle {   \sum_{j=1}^{d}M^{21}_{j}\partial _j U^{I0} =
   Q_ \nu(U^{I0},0)U^{II0}+D^{II}(U^{I0}, 0)},
  \end{array}
\right.
\label{24}
\end{equation}
in the sense of distribution.
\label{t5}
\end{theorem}
\begin{proof}
(\ref{21}), (\ref{22}) follow from (i) of the Theorem \ref{t4} and
the energy estimate implies that
$$
\begin{pmatrix}
 U^{I}_{t}+\displaystyle{\sum_{j=1}^{d}M^{12}_{j}\partial _j U^{II}} \\
\displaystyle {\sum_{j=1}^{d}M^{21}_{j}\partial _j U^{I}}
\end{pmatrix} \quad \quad \text{is relatively compact in $(H^{-1}_{loc})^2$}.
$$
In order to fit into the framework of Theorem \ref{cc} we set
$l^\varepsilon =(U^I, U^{II})$, then the characteristic manifold $\Lambda$
is given by
\begin{displaymath}
    \Lambda =\left\{(\lambda, \mu)\in \ERRE ^k \times \ERRE ^{N-k}
        ~|~ \exists
    \xi\in
    \mathbb{R}^{d+1}\setminus\left\{0\right\},
    B(\xi,\lambda)=0\right\}
\end{displaymath}
where
\begin{displaymath}
    B(\xi,\lambda)=\begin{pmatrix}
\xi_0\lambda+
\left(\displaystyle{\sum _{j=1}^d M^{22}_j \xi _j}\right) \mu\\
\left(\displaystyle{\sum _{j=1}^d M^{21}_j (x)\xi _j}\right) \lambda
\end{pmatrix}.
\end{displaymath}
By using the hypothesis (B.4) and by defining
 $\Theta=\frac 1 2\left(\begin{smallmatrix}
    I_{k \times k} & 0 \\
     0            & 0  \\
    \end{smallmatrix}\right) $,
 we have $f(\lambda)={\lambda}^{T}\Theta\lambda$ and, of
course, $f_{|\Lambda}\equiv 0$.
We apply the Theorem \ref{cc} to show
$$
(U^I)^2 \rightharpoonup  (U^{I0})^2 \qquad \text{ in the sense of measure }
$$
and finally
$$
U^I \longrightarrow U^{I0} \qquad \text{strongly in
$L^{2}_{loc}(\ERRE ^d \times \ERRE_+)$}.
$$
So we can pass to the limit into the nonlinear terms of the system (\ref{18}).
\end{proof}
\begin{remark}
An equivalent form of system (\ref{24}) is given by (setting $U^{I0}=U$)
\begin{equation}
\begin{split}
 U_{t}&+\sum_{j=1}^{d}M^{12}_{j}\partial _j\left(Q_ \nu^{-1}(x,U,0)
\sum_{k=1}^{d} M^{21}_{k}\partial _{k} U \right)\\&=
\sum_{j=1}^{d}M^{12}_{j}\partial_j(Q_ \nu^{-1}(x,U,0) D^{II}(U,0))
\\
 & +D^I_ \nu(U^{I},0)Q_ \nu^{-1}(x,U,0)
\left[\sum_{k=1}^{d}M^{21}_{k}\partial _k U - D^{II}(U, 0)\right].
\end{split}
\label{25}
\end{equation}
The proof of the parabolicity of (\ref{25}) follows the same arguments used in
the variable coefficients case.
\end{remark}
\end{subsection}

\end{section}

\begin{section}{Approximation of given parabolic systems}
In this last section we want to apply the theory of the previous
ones to approximate a generic given parabolic system, provided the 
Petrowski condition hold. 
In fact here we reconstruct a parabolic system by means of a
suitable larger semilinear hyperbolic system that relaxes on it.
Two very large classes of parabolic systems will be taken in 
consideration. The former are quasilinear parabolic systems in divergence 
form and the latter are the so called ``Reaction- diffusion'' systems.

\begin{subsection}{Quasilinear case}

We consider now the following quasilinear system in divergence form
\begin{equation}
U_t +\sum _{i=1}^{d}\partial _{i} \left( F_{i}(U)-\sum _{j=1}^{d}B_{ij}(U)\partial 
_{j}U\right )=G(U) \label{70}
\end{equation}
where $x \in \ERRE ^d$, $t \in \ERRE_+$, $U=U(x,t) \in \ERRE ^k$,
$N=(d+1)k$ ( then $k\leq\frac{N}{2} )$. Let us denote by $\F(U)\in 
\mathcal{M}_{k\times d}$, $(\F(U))_{i}=F_{i}(U)$ for any $i=1,\ldots, d$ 
and $\B(U)\in 
\mathcal{M}_{kd\times kd}$, $(\B(U))_{ij}=B_{ij}(U)$ for any $i,j=1,\ldots, 
d$, then we assume.
\begin{itemize}
\item[\bf(C.1)]
$B_{ij}(\cdot)\in C^1 (\ERRE ^k;\mathcal{M}_{k\times d})$, for any  $i,j=1,\ldots, d$
and $\sum_{i,j}B_{ij}(U) \lambda _i\lambda _j \geq c_0 |\lambda|^2I $ for any 
$\lambda \in \ERRE ^d$ (strong parabolicity), moreover
$\B^{-1}(U)$ is bounded on $U$,
\item[\bf(C.2)] $F_{i}(\cdot)\in C^1 (\ERRE ^k;\mathcal{M}_{k\times d})$, 
for any $i=1,\ldots, d$, $\F(0)=0$ 
and $\B^{-1}(U)\F(U)$ is lipschitz   on $U$.
\item[\bf(C.3)] $G(\cdot)\in Lip(\ERRE ^k; \ERRE ^k)$.
\end{itemize}
We have the following theorem
\begin{theorem}
Let us consider the system (\ref{70}), let us suppose  hypotheses
(C.1), (C.2), (C.3) holds, then the solutions of the system
\begin{equation}
\left\{\begin{array}{ll} 
Z^I_t +\displaystyle{\frac{1}{\varepsilon}\dive \Z^{II}}= G(Z^{I})\\ 
\Z^{II}_t+
\displaystyle{\frac{1}{\varepsilon}DZ^{I}}= \displaystyle{-\frac{1}{\varepsilon
^2}}\B^{-1}(Z^I)\Z^{II}+\frac{1}{\varepsilon}\B^{-1}(Z^I)\F(Z^I),
\end{array}
\right. \label{71}
\end{equation}
where   $(x,t)\in \ERRE ^d\times \ERRE_+$,
$Z^I=Z^I(x,t)\in\ERRE ^k$, $\Z^{II}=\Z^{II}(x,t)\in \mathcal{M}_{k\times d}$,
approximate the system (\ref {70}) in the sense of Theorem \ref{t2} and the Theorem
\ref{t3}.
\end{theorem}
\begin{proof}
It can be easily shown  that  (\ref{71}) is an hyperbolic symmetric system. 
By rescaling the variables,
as in (\ref{30}), the  system (\ref{71}) transforms into
\begin{equation}
\left\{\begin{array}{ll} U^I_t+
\dive \U^{II}=G(U^{I})
\\ \varepsilon ^2 \U^{II}_t+
DU^{I}=
-\B^{-1}(U^I)\U^{II}+\B^{-1}(U^I)\F(U^I).
\end{array}
\right. \label{72}
\end{equation}
The conditions (B.1), (B.2) of the 
Section 4 are satisfied by the system (\ref{72}).Now by 
setting $D^{II}(Z^I,\Z^{II}\!)=
\B^{-1}(Z^I)\F(Z^{I})$ we can easily verify (B.3) and by
denoting  $Q(x,Z^I,\Z^{II})=-\B^{-1}(Z^I)\Z^{II}$ we get obviously the
condition (D). We can apply the Theorem \ref{t3} and we obtain
that the solutions to (\ref{72}) satisfy, as 
$\varepsilon \downarrow 0$ (setting $U^{I0}=U$) 
$$U_t+\dive (\F(U)+\B(U)DU)=G(U).$$     
\end{proof}
\end{subsection}

\begin{subsection}{Reaction - Diffusion type systems}
Let us consider the following reaction-diffusion type system
\begin{equation}
U_t=\sum_{j,k=1}^{d}A_{j,k}(x)\partial _j\partial _k U +f(U)\label{50}
\end{equation}
where $x \in \ERRE ^d$, $t \in \ERRE_+$, $U=U(x,t) \in \ERRE ^k$. 
We make the following hypotheses.
\begin{itemize}
\item[\bf(D.1)]
$A_{j,k}(\cdot) \in C^\infty (\ERRE ^d ;\mathcal{M}_{k\times k})$, for
any  $j,k=1,\ldots d$ and $\sum _{j,k=1}^{d}A_{j,k}(x) \lambda
_j\lambda _k \geq c_0 |\lambda|^2I $, for any $\lambda \in \ERRE
^d$ (strong parabolicity).
\item[\bf(D.2)] $f(\cdot)\in Lip(\ERRE ^k; \ERRE ^k)$.
\end{itemize}
In order to approximate system \eqref{50}, we define the following linear 
operator;\\ $B_{j}(x):\ERRE^{kd}\longrightarrow \ERRE^{kd}$, 
$$B_{j}(x)(W)=\sum_{k=1}^{d}A_{j,k}(x)W_{k}\qquad \text{for any $j=1,\ldots 
,d$}$$
Now system \eqref{50} can be written in the equivalent form
$$U_t=\sum_{j=1}^{d}B_{j}(\partial _{j}DU)+f(U).$$
We have the following theorem.
\begin{theorem}
Let us consider the system (\ref{50}), let us suppose the hypotheses
(D.1), (D.2) hold, then the solutions of the system
\begin{equation}
\left\{\begin{array}{ll}
Z^I_t
+\displaystyle{\frac{1}{\varepsilon}}\sum_{j=1}^{d}B_{j}\partial 
_{j}\Z^{II}=f(Z^{I})\\
 \Z^{II}_t+\displaystyle{\frac{1}{\varepsilon}}\sum_{j=1}^{d}B_{j}^{T}
\partial _{j}Z^{I}=\displaystyle{-\frac{1}{\varepsilon ^2}}
\A (x)\Z^{II}
      \end{array}
\right. 
\label{531}
\end{equation}
where  $(x,t)\in \ERRE ^d\times \ERRE_+$, $Z^I=Z^I(x,t)\in\ERRE
^k$, $\Z^{II}=\Z^{II}(x,t)\in\ERRE ^{kd}$, $\A(x)\in \mathcal{M}_{kd\times 
kd}$, $(A(x))_{j,k}=A_{j,k}(x)$ approximate the system
(\ref {50}) in the sense of the Theorem \ref{t2} and the Theorem
\ref{t3}.
\label{t531}
\end{theorem}
\begin{proof}
    The system \eqref{531} is symmetric and hyperbolic. By rescaling the 
    variables as in (\ref{30}) it transforms into
 \begin{equation}
\left\{\begin{array}{ll}
U^I_t
+\displaystyle{\sum_{j=1}^{d}B_{j}\partial 
_{j}\U^{II}}=f(U^{I})\\
 \varepsilon ^{2}\U^{II}_t+\displaystyle{\sum_{j=1}^{d}B_{j}^{T}
\partial _{j}U^{I}}=-\A (x)\U^{II}
      \end{array}
\right. 
\label{532}
\end{equation} 
By setting $Q(x,Z^I, \Z^{II})=-\A(x)\Z^{II}$, the dissipativity 
conditions (D) together with the hypotheses (B.1) and 
(B.2) of Section 4 are immediately satisfied. Now, we can apply the
Theorems \ref{t2} and \ref{t3}, then since 
$\displaystyle{\sum_{j=1}^{d}B_{j}^{T}
\partial _{j}U^{I}}=\A(x)DU^{I}$, as $\varepsilon \downarrow 0$, 
the solution of (\ref{532})
satisfy the system \eqref{50}.
    \end{proof}
\begin{remark} 
The Theorem \ref{t531} can be applied,  with slight modifications,
to the more general case of Petrowski 
parabolic systems in the sense recalled in the Section 2.
In fact let us consider the following system
\begin{equation}
U_t=\sum_{j,k=1}^{d}C_{j,k}(x)\partial _j\partial _k U +G(U)\label{533}
\end{equation}
where $x \in \ERRE ^d$, $t \in \ERRE_+$, $U=U(x,t) \in \ERRE ^k$ and the 
matrices $C_{j,k}$ satisfy the Petrowski parabolicity condition
given in the Section 2.
Let us consider the matrix $P_{0}$ from the previous definition, denote by
$$W=P_{0}U$$
then, the system \eqref{533} transforms into a strongly parabolic system, so we 
can apply the previous Theorem \ref{t531}.
\end{remark}
To approximate the system \eqref{50} we can also follow a different 
approach by using the pseudodifferential theory. In fact, let us denote by
\begin{equation}\label{51}
A(x,\xi)=-\sum _{j,k}A_{j,k}(x)\xi _j\xi _k,
\end{equation}
$\xi \in \ERRE ^d$, the principal symbol of (\ref{50}) and by
$A(x,D)=OPA(x,\xi)$. Now we set
\begin{equation}\label{52}
B(x,D)=OP[(-A(x,\xi))^{1/2}]
\end{equation}
We have the following theorem.
\begin{theorem}
Let us consider the system (\ref{50}), suppose that hypotheses
(D.1), (D.2) hold, then the solutions of the system
\begin{equation}
\left\{\begin{array}{ll}
Z^I_t
+\displaystyle{\frac{1}{\varepsilon}}B(x,D)Z^{II}=f(Z^{I})\\
 Z^{II}_t-\displaystyle{\frac{1}{\varepsilon}}B(x,D)Z^{I}=
\displaystyle{-\frac{1}{\varepsilon ^2}}Z^{II},
      \end{array}
\right. \label{53}
\end{equation}
where  $(x,t)\in \ERRE ^d\times \ERRE_+$, $Z^I=Z^I(x,t)\in\ERRE
^k$, $Z^{II}=Z^{II}(x,t)\in\ERRE ^{k}$, $N=2k$ approximate the system
(\ref {50}) in the sense of the Theorem \ref{t2} and the Theorem
\ref{t3}.
\end{theorem}
\begin{proof}
By rescaling the variables as in (\ref{30}) the system (\ref{50})
transforms into
\begin{equation}
\left\{\begin{array}{ll} U^I_t +B(x,D)U^{II}=f(U^{I})\\ \\ \varepsilon ^2
U^{II}_t-B(x,D)U^{I}=- U^{II}
      \end{array}
\right. \label{54}
\end{equation}
By \ using \ (D.1) \ we \ have   \ (\ref{54})\  is  \ an 
hyperbolic system. Let\ 
us \  denote  \ by   $M^{12}(x,D)=B(x,D)$, $M^{21}(x,D)=-B(x,D)$, by using 
the hypothesis(C.1), the conditions
(B.1) and (B.2) are satisfied. It can also be easily
verified that the symmetrizer of (\ref{53}) is given by the matrix
$I_{N\times N}$. Finally we have to verify the dissipativity
condition (D). Setting $Q(x,Z^I, Z^{II})=-Z^{II}$ we get the hypothesis (D)
is satisfied. Now we can apply Theorem \ref{t2} and
Theorem \ref{t3} and we obtain that the solution of (\ref{54})
satisfy as $\varepsilon \downarrow 0$ (by setting $U^{I0}=U$) $$U_t
-A(x,D)U=0.$$ 
\end{proof}
\end{subsection}

\end{section}
\bibliographystyle{amsplain}

\end{document}